\documentclass{article}
\usepackage{amsmath,amssymb,amsthm,amscd, mathtools, latexsym}
\usepackage{enumerate,varioref, dsfont, fancyhdr}
\usepackage{tikz-cd}
\usepackage{enumitem}
\usepackage{multicol}
\usepackage{hyperref}
\usepackage{url}
\usepackage{textcomp}

\newtheorem{Thm}{Theorem}[section]
\newtheorem{Lem}[Thm]{Lemma}
\newtheorem*{Lem*}{Lemma}
\newtheorem{Prop}[Thm]{Proposition}

\newtheorem{Cor}[Thm]{Corollary}
\newtheorem{Conj}[Thm]{Conjecture}
\newtheorem{Rem}[Thm]{Remark}

\newtheorem{Def}[Thm]{Definition}

{\theoremstyle{plain}

\newtheorem*{Ass*}{Assumption}

\newtheorem{Not}[Thm]{Notation}

}

\def\zit{{\mathbb Z}}

\def\pit{{\mathbb P}}

\def\0{{\mathcal O}}

\def\End{\mathop{\rm End}\nolimits}

\def\H{{\mathcal H}}

\def\C{{\mathcal C}}

\def\M{{\mathcal M}}

\usepackage[title]{appendix}

\usepackage[left=3.5cm, right=3.5cm]{geometry}

\begin{document}
\title{Rost nilpotence and higher unramified cohomology}
\author{H. Anthony Diaz}
\newcommand{\Addresses}{{\bigskip \footnotesize
\textsc{Department of Mathematics, Washington University, St. Louis, MO 63130} \par \nopagebreak
\textit{Email address}: \ \texttt{humberto@wustl.edu}}}

\date{}
\maketitle

\begin{abstract}
\noindent We develop an approach to proving the Rost nilpotence principle involving higher unramified cohomology. Along with a decomposition-of-the-diagonal technique, we use this to prove the principle for certain threefolds over a perfect field.
 
\end{abstract}
\date{}

\section*{Introduction}

\noindent Let $X$ be a smooth projective scheme over a field $k$. Also, let $\M_{k}$ be the category of Chow motives over $k$ and denote the corresponding Chow motive of $X$ by $M(X)$. The {\em Rost nilpotence principle} predicts that for any field extension $E/k$, the kernel of the extension of scalars map
\begin{equation} \End_{\M_{k}} (M(X)) \to \End_{\M_{E}} (M(X_{E}))\label{extend} \end{equation}
consists of nilpotent correspondences. Rost first proved this for a smooth projective quadric over a field of characteristic $\neq 2$ in \cite{R} (see also \cite{Br}). A consequence is that the Chow motive of a smooth quadric can be decomposed as a direct sum of (twisted) motives of anisotropic quadrics, and this played an important role in the proof of the Milnor conjecture by Voevodsky \cite{V}.\\
\indent The Rost nilpotence principle is desirable, as it would imply for instance that Chow motives do not vanish upon passage to field extensions (in this sense, it may be viewed as a torsion analogue of a well-known nilpotence conjecture for rational correspondences \cite{Ki}). It is conjectured to hold in general and has been proved in several other important cases. Chernousov, Gille and Merkurjev \cite{CGM} proved that it holds for projective homogeneous varieties. Moreover, using Rost cycle modules \cite{R1}, Gille showed that the Rost nilpotence principle holds for geometrically rational surfaces over perfect fields \cite{Gi1} and for smooth, projective, geometrically integral surfaces over fields of $\text{char } k = 0$ \cite{Gi}. (Using different methods, he also proved it for smooth projective, geometrically integral threefolds which are birationally isomorphic to toric models \cite{Gi2}.)\\
\indent Another approach to proving Rost nilpotence was developed by Rosenschon and Sawant \cite{RS}. Their approach involves \'etale motivic cohomology groups $H^{m}_{L} (X, \zit(n))$ which in many ways are better behaved than the usual motivic cohomology groups (see, for instance, \cite{RSr}). In particular, for any finite Galois extension $E/k$ there is a Hochschild-Serre spectral sequence:
\[ E_{2}^{p,q}=H^{p} (Gal(E/k), H^{q}_{L} (X_{E}, \zit(n)) \Rightarrow H^{p+q}_{L} (X, \zit(n)) \]
Moreover, using the triangulated category of \'etale mixed motives, this spectral sequence is functorial for the action of correspondences after inverting the exponential characteristic of $k$ (see \cite{RS} \S3). Using this, they are able to prove that if $\gamma \in CH^{d} (X \times X)$ lies in the kernel of (\ref{extend}), then the action of $\gamma$ on the \'etale motivic cohomology groups of $X$ (and its products) is nilpotent. As a consequence, they obtain the Rost nilpotence principle for smooth projective surfaces in characteristic $0$ (and with some more work, birationally ruled threefolds in characteristic $0$).\\
\indent Our goal in this paper will be to explore the relationship between Rost nilpotence and higher unramified cohomology. The usual unramified cohomology groups $H^{*}_{nr} (X)$ are well-known and since \cite{CTO} have been quite useful in proving that certain varieties are not rational (or stably rational). The right derived analogues of these groups are the higher unramified cohomology groups.
The na\"ive view is that the higher unramified cohomology groups control the extent to which motivic cohomology and \'etale motivic cohomology fail to coincide. Thus, given the main result of \cite{RS}, one would expect the Rost nilpotence principle to follow from a statement about correspondences acting on unramified cohomology. In this direction, our main result (Theorem \ref{main}) shows that a correspondence $\gamma$ that lies in the kernel of (\ref{extend}) is nilpotent, provided that the action of $\gamma$ on certain higher unramified cohomology groups is nilpotent. In particular, since the assumption of Theorem \ref{main} is automatically satisfied for surfaces, we are able to extend the Rost nilpotence principle to surfaces over a perfect field. Then, using a Bloch-Srinivas-type argument, we are able to prove the following as a consequence:
\begin{Thm}[=Corollary \ref{three}] Suppose that $X$ is a smooth projective scheme of dimension $\leq 3$ over a perfect field $k$ whose Chow group $CH_{0}$ is universally supported in dimension $\leq 2$ (in the sense of \cite{BS}; see also Definition \ref{def}). Then, $X$ satisfies the Rost nilpotence principle. 
\end{Thm}

\noindent Since the Rost nilpotence principle is known for surfaces over a perfect field (after inverting the exponential characteristic) by the main result of \cite{RS}, we focus on the case that $k$ has positive characteristic and fill in the gap involving $p$-primary torsion. For this, we make use of the higher unramified cohomology groups for logarithmic Hodge-Witt cohomology, as well as that of the Hochschild-Serre spectral sequence for \'etale cohomology groups. The idea will be to prove that any cycle in the kernel of (\ref{extend}) is nilpotent in \'etale cohomology and then to use the local-to-global spectral sequence to obtain a nilpotence statement in the Chow group. The local-to-global spectral sequence in this case reduces to a long exact sequence by some non-trivial results of Gros and Suwa (discussed below). Since the cycle class map is not guaranteed to be injective even in low degree, we make use of a Bockstein-type sequence that overcomes this difficulty. The proofs of the applications that follow Theorem \ref{main} are ``decomposition-of-the-diagonal" arguments in the spirit of \cite{BS}.

\subsection*{Acknowledgements}
The author would like to thank Bruno Kahn for taking the time to read several drafts of an earlier paper that inspired this one; his comments were indispensable. The author also thanks Anand Sawant for his interest and Fr\'ed\'eric D\'eglise for a clarification. Finally, the author thanks the referees for discovering some gaps in the previous version and for correctly identifying how to fix them, as well as for making numerous other suggestions.
\subsection*{Conventions}
Throughout, we will assume that all schemes are separated and of finite type over a field $k$. When we say that a scheme $X$ has dimension $d$, we mean that it is equi-dimensional and that all the irreducible components have dimension $d$. In particular, a {\em surface} will mean a scheme of dimension $2$; similarly, a {\em threefold} will mean a scheme of dimension $3$. In the case that $X$ is irreducible, we let $k(X)$ denote the function field and $\eta$ the generic point. For any field extension $L/k$ and $X$ a scheme over $k$, we write $X_{L} = X \times_{k} L$. We also let $\overline{k}$ denote the algebraic closure and $\overline{X} = X \times_{k} \overline{k}$.\\
\indent Chow groups (indexed by codimension) will be denoted by $CH^{i}(-)$ and have integral coefficients. In the final section, we also index by dimension. If $X$ and $Y$ are both smooth and projective schemes; let $Y_{1}, \ldots Y_{m}$ be the irreducible components of $Y$ with respective dimensions $d_{1}, \ldots d_{m}$. Then, we denote by
\begin{equation} Cor^{i} (X, Y) := \bigoplus_{k=1}^{m} CH^{i+d_{k}} (X \times Y_{k}) \label{cor-def} \end{equation}
the group of correspondences of degree $i$ and $Cor^{*} (X, Y) = \bigoplus_{i \in \zit} Cor^{i} (X \times Y)$ the graded ring of correspondences, where muliplication is given by the composition operation $\circ$ defined in \cite{F} Ch. 16.\\
\indent The notation $\ell$ will be reserved for a prime number. For an Abelian group $A$ and $n \geq 1$ an integer, $A[n]$ denotes the $n$-torsion of $A$; in particular, for a prime number $\ell$, $A[\ell^{\infty}]$ denotes the $\ell$-primary torsion of $A$.\\

\section{Preliminaries}
\noindent Throughout this section, $X$ will denote a smooth scheme over a perfect field $k$ of characteristic $p>0$ (unless otherwise specified).
\subsection{Logarithmic Hodge-Witt cohomology}
\noindent We recall some important facts about logarithmic Hodge-Witt cohomology. For $r \geq 0$, let $W_{r}\Omega_{X}^{*}$ be the de Rham-Witt complex defined by Illusie \cite{I} Ch. 1; this is a complex of Zariski sheaves over $X$. Also, the logarithmic Hodge-Witt complex over $X$, $W_{r}\Omega_{X, log}^{*}$ is the complex of Zariski sheaves over $X$ defined in op. cit. Ch. 1 \S 5.7. We will denote the \'etale sheafification by $\nu_{r}(n)_{X}:= \alpha^{*}W_{r}\Omega_{X, log}^{n}$, where $\alpha$ denotes the forgetful functor from the \'etale site to the Zariski site. This sheaf is locally generated by differentials of the form 
\[\frac{d\tilde{x}_{1}}{\tilde{x}_{1}} \wedge \ldots \wedge \frac{d\tilde{x}_{n}}{\tilde{x}_{n}} \]
where $\tilde{x}_{i}$ denotes a Teichm\"uller representative of $x_{i} \in \mathcal{O}_{X}^{*}$. 
Note also that by (5.1.1) of op. cit. there is an isomorphism of cohomology groups:
\[ H^{*}_{Zar} (X, W_{r}\Omega_{X, log}^{n}) \cong  H^{*}_{\text{\'et}} (X, \nu_{r}(n)_{X})\]
As is standard (see, for instance, \cite{GL}), we shift by $n$ and consider the complex of \'etale sheaves over $X$ concentrated in degree $n$:
\[ (\zit/p^{r})_{X}(n) := \nu_{r}(n)_{X}[-n] \]
With this definition, we have $H^{i}_{\text{\'et}} (X, (\zit/p^{r})_{X}(n)) = H^{i-n}_{\text{\'et}} (X, \nu_{r}(n)_{X})$. Additionally, note that when $n=0$ we recover the usual locally-constant \'etale sheaf $\zit/p^{r}$ (see p. 597 of \cite{I}).\\
 \indent For any morphism of smooth schemes over $k$ $f: Y \to X$, there is a canonical map in the derived category of \'etale sheaves of Abelian groups over $X$:
 \begin{equation} (\zit/p^{r})_{X}(n) \to Rf_{*}(\zit/p^{r})_{Y}(n)\label{pull}\end{equation}
(see \cite{G} (1.2.1)); applying \'etale hypercohomology to (\ref{pull}) yields the {\em pull-back} map
\[ H^{i}_{\text{\'et}} (X, (\zit/p^{r})_{X}(n)) \xrightarrow{f^{*}} H^{i}_{\text{\'et}} (Y, (\zit/p^{r})_{Y}(n)) \]
From (1.2.3) of op. cit., there is also an associative and anti-commutative product map 
\[ \zit/p^{r}(n) \otimes \zit/p^{r}(n') \to \zit/p^{r}(n+n') \] 
(induced by the product on $W_{r}\Omega_{X}^{*}$); the induced map on cohomology
 \begin{equation} H^{i}_{\text{\'et}} (X, (\zit/p^{r})_{X}(n)) \otimes H^{i'}_{\text{\'et}} (X, (\zit/p^{r})_{X}(n')) \xrightarrow{\cup} H^{i+i'}_{\text{\'et}} (X, (\zit/p^{r})_{X}(n+n'))\label{cup-prod} \end{equation}
is the {\em cup product}. Finally, for proper morphisms $f: Y \to X$ of smooth schemes over $k$ of relative dimension $d$, there is the trace map (see p. 5 of op. cit.):
 \begin{equation} Rf_{*}(\zit/p^{r})_{Y}(n) \to (\zit/p^{r})_{X}(n-d)[-2d] \label{push}\end{equation}
induced by the trace map for the corresponding de Rham-Witt complexes. This induces the {\em push-forward} map on cohomology
 \[ H^{i}_{\text{\'et}} (Y, (\zit/p^{r})_{Y}(n)) \xrightarrow{f_{*}} H^{i-2d}_{\text{\'et}} (X, (\zit/p^{r})_{X}(n-d)) \]
\indent The above operations of pull-back, push-forward and cup product allow one to view cohomology classes as {\em cohomological correspondences} that may be composed. Indeed, suppose that $X$, $Y$ and $Z$ are smooth projective schemes of dimension $d_{Y}$ and $d_{Z}$ over $k$ and suppose that 
\[ \Gamma \in H^{i+2d_{Y}}_{\text{\'et}} (X\times Y, (\zit/p^{r})_{X \times Y}(n+d_{Y})), \ \Gamma' \in H^{j+2d_{Z}}_{\text{\'et}} (Y\times Z, (\zit/p^{r})_{Y \times Z}(s+d_{Z})) \]
Then, we define their composition in the usual way (see, for instance, \cite{F} p. 305-307):
\begin{equation} \Gamma' \circ \Gamma := \pi_{XZ*}(\pi_{XY}^{*}\Gamma\cdot\pi_{YZ}^{*}\Gamma') \in  H^{i+j+2d_{Z}}_{\text{\'et}} (X\times Z, (\zit/p^{r})_{X \times Z}(n+s+d_{Z}))\end{equation}
One checks that $\circ$ is associative and, hence, induces a ring structure on logarithmic Hodge-Witt cohomology. Indeed, as in Prop. 16.1.1 of loc. cit., this follows formally from the functoriality of the pull-back, push-forward and the projection formula, all of which are established in \cite{G}. Given a cohomological correspondence $\Gamma$ as above, there is an induced action on logarithmic Hodge-Witt cohomology given by the usual formula:
\[ \Gamma_{*}(\alpha):= \pi_{Y*}(\pi_{X}^{*}\alpha\cdot \Gamma) \in H^{*+i} (Y, (\zit/p^{r})_{Y}(m+n)) \]
for $\alpha \in H^{*} (X, (\zit/p^{r})_{X}(m))$. The associativity of $\circ$ immediately implies that $(\Gamma'\circ\Gamma)_{*} = \Gamma'_{*}\circ\Gamma_{*}$.

\subsection{Hochschild-Serre spectral sequence} 

\noindent For any bounded complex of \'etale sheaves $C^{*}$ over $X$, we let $\mathbb{H}^{*}_{\text{\'et}} (X, C^{*})$ denote the corresponding \'etale hypercohomology group (\cite{M} Appendix C). Then, there is the Hochschild-Serre spectral sequence:
\begin{equation} E_{2, C^{*}}^{s,t} = H^{s} (k, \mathbb{H}_{\text{\'et}}^{t} (\overline{X}, C^{*})) \Rightarrow \mathbb{H}^{s+t}_{\text{\'et}} (X, C^{*}) \label{HS}\end{equation}
The convergence of this spectral sequence follows from the fact that $C^{*}$ is bounded (see, for instance, \cite{K} 2C). This sequence is functorial in the sense that if $C^{*} \to D^{*}$ is a map of bounded complexes of \'etale sheaves over $X$, then there is an induced map of the spectral sequences (in the sense of \cite{M} p. 307-308) for which the map on abutments is the hypercohomology map:
\[ \mathbb{H}^{s+t}_{\text{\'et}} (X, C^{*}) \to \mathbb{H}^{s+t}_{\text{\'et}} (X, D^{*}) \]
This map of spectral sequences induces a map on the descending filtration of the abutment.\\ 
\indent We will consider (\ref{HS}) in the case that $C^{*} = \zit/p^{r}(n)$. Then, we have the Hochschild-Serre spectral sequence:
\[ H^{s} (k, H_{\text{\'et}}^{t} (\overline{X}, \zit/p^{r}(n))) \Rightarrow H^{s+t}_{\text{\'et}} (X, \zit/p^{r}(n)) \]
Also, letting $F^{*}$ denote the corresponding filtration, we have:
\[F^{0}H^{m}_{\text{\'et}} (X, \zit/p^{r}(n)) = H^{m}_{\text{\'et}} (X, \zit/p^{r}(n)), \ F^{m+1}H^{m}_{\text{\'et}} (X, \zit/p^{r}(n)) = 0 \]
for degree reasons. Moreover, there is a natural short exact sequence:
\begin{equation} 0 \to F^{1}H^{m}_{\text{\'et}} (X, \zit/p^{r}(n)) \to H^{m}_{\text{\'et}} (X,\zit/p^{r}(n)) \to H^{m}_{\text{\'et}} (\overline{X}, \zit/p^{r}(n))^{G_{k}}\label{other-edge} \end{equation}
where $G_{k} = Gal(\overline{k}/k)$ denotes the absolute Galois group of $k$. 
\begin{Prop}\label{Galois-props} Suppose that $X$ is a smooth and projective scheme over $k$. Then, the filtration $F^{*}$ is functorial with respect to pull-backs, pushforwards and cup product; i.e., for any morphism $f: Y \to X$ (of relative dimension $d$), we have:
\begin{enumerate}[label=(\alph*)]
\item $f^{*}(F^{i}H^{m}_{\text{\'et}} (X, \zit/p^{r}(n)))  \subset F^{i}H^{m}_{\text{\'et}} (Y, \zit/p^{r}(n))$;
\item $f_{*}(F^{i}H^{m}_{\text{\'et}} (Y, \zit/p^{r}(n)))  \subset F^{i}H^{m-2d}_{\text{\'et}} (X, \zit/p^{r}(n-d))$;
\item $F^{i}H^{m}_{\text{\'et}} (X, \zit/p^{r}(n))\cup F^{j}H^{s}_{\text{\'et}} (X, \zit/p^{r}(t)) \subset F^{i+j}H^{m+s}_{\text{\'et}} (X, \zit/p^{r}(n+t))$
\end{enumerate}
\begin{proof} As noted previously, the pull-back and push-forward maps are defined on the level of the derived category of complexes of bounded \'etale sheaves:
\[ (\zit/p^{r})_{X}(n) \to Rf_{*}(\zit/p^{r})_{Y}(n), \ Rf_{*}(\zit/p^{r})_{Y}(n) \to (\zit/p^{r})_{X}(n-d)[-2d]\] 
which implies by the above paragraphs that the corresponding pull-back and push-forward maps on logarithmic Hodge-Witt cohomology respect the filtration. For the compatibility with cup product, note that the Hochschild-Serre spectral sequence {\em has cup products} in the terminology of \cite{Sw} p. 184; this fact can be deduced (for instance) from the argument given for Cor. 3.11 of op. cit.
\end{proof}
\end{Prop}
\noindent As a consequence, we have the following result about the vanishing of sufficiently high powers of certain cohomological correspondences. We will need this in the proof of our main result.
\begin{Cor}\label{cor-van} For $X$ a smooth projective scheme of dimension $d$ over $k$, let 
\[ \Gamma \in \text{ker } \{ H^{2d}_{\text{\'et}} (X \times X, \zit/p^{r}(d)) \to H^{2d}_{\text{\'et}} (\overline{X} \times \overline{X}, \zit/p^{r}(d)) \} \] 
be a cohomological correspondence. Then, $\Gamma^{\circ m} = 0 \in H^{2d}_{\text{\'et}} (X \times X, \zit/p^{r}(d))$ for $m \geq d+1$.
\begin{proof} From the previous result, it follows that $\Gamma^{\circ m} \in F^{m}H^{2d}_{\text{\'et}} (X \times X, \zit/p^{r}(n))$. Note that $E_{2}^{m,2d-m} = 0$ unless $m \leq d$, since for $m > d$, we have
\[ H^{2d-m}_{\text{\'et}} (X \times X, \zit/p^{r}(d)) = H^{d-m}_{\text{\'et}} (X \times X, \nu_{r}(d)) = 0 \]
So, it follows that $F^{m}H^{2d}_{\text{\'et}} (X \times X, \zit/p^{r}(n)) =0$ for $m \geq d+1$. This gives the desired result.
\end{proof}
\end{Cor}
\begin{Rem} If $k$ is a field of cohomological dimension $r$, then we actually have $\Gamma^{\circ m} = 0$ for $m \geq r+1$.
\end{Rem}
\subsection{Bloch-Ogus spectral sequence}

\indent Now, we let $\mathcal{H}_{X}^{q}(\zit/p^{r}(n))$ denote the Zariski sheaf associated to the presheaf of Abelian groups over $X$ given by $U \mapsto H^{q}_{\text{\'et}} (U, \zit/p^{r}(n))$. (From this point forward, we will drop the $X$ in the subscript and use the notation $\zit/p^{r}(n)$ and $\nu_{r}(n)$ instead of $(\zit/p^{r})_{X}(n)$ and $\nu_{r}(n)_{X}$, when there is no danger of ambiguity.) By general properties of right-derived functors (for instance, \cite{M} Ch. 3 Prop. 1.13), we have the identification:
\begin{equation} \mathcal{H}_{X}^{q}(\zit/p^{r}(n)) \cong R^{q}\alpha_{*}\zit/p^{r}(n) = R^{q-n}\alpha_{*}\nu_{r}(n)\label{formal}\end{equation} 
noting the above-mentioned degree shift. The corresponding Zariski cohomology groups:
\[H^{i}_{Zar} (X,\mathcal{H}_{X}^{q}(\zit/p^{r}(n))) \]
are what we refer to as the {\em higher unramified cohomology groups}. Then, there is the corresponding Grothendieck spectral sequence abutting to logarithmic Hodge-Witt cohomology:
\begin{equation} E_{2}^{i,q} = H^{i}_{Zar} (X,\mathcal{H}_{X}^{q}(\zit/p^{r}(n))) \Rightarrow H^{i+q}_{\text{\'et}} (X, (\zit/p^{r})_{X}(n))\label{pre-coh}  \end{equation}
known as the {\em Bloch-Ogus spectral sequence}. By \cite{GS} item 1.10, (\ref{formal}) vanishes for $q \neq n,n+1$, so $E_{2}^{i,q} = 0$ unless $q = n,n+1$. Consequently, the Bloch-Ogus spectral sequence reduces to a long exact sequence (\cite{GS} (1.14)):
\begin{equation}  \begin{split}\cdots\to & H^{m}_{Zar} (X, \mathcal{H}^{n}_{X}(\zit/p^{r}(n)))
\xrightarrow{e^{m, n}} H^{m+n}_{\text{\'et}} (X, \zit/p^{r}(n)) \to H^{m-1}_{Zar} (X, \mathcal{H}^{n+1}_{X}(\zit/p^{r}(n)))\\ \xrightarrow{d_{2}^{m-1,n+1}} & H^{m+1}_{Zar} (X, \mathcal{H}^{n}_{X}(\zit/p^{r}(n))) \xrightarrow{e^{m+1, n}} H^{m+n+1}_{\text{\'et}} (X, \zit/p^{r}(n)) \to \cdots \label{long}
\end{split}\end{equation}
\indent Finally, we observe that there is a natural product structure:
\begin{equation} \mathcal{H}^{m}_{X}(\zit/p^{r}(n)) \otimes \mathcal{H}^{q}_{X}(\zit/p^{r}(s)) \xrightarrow{\cup} \mathcal{H}^{m+q}_{X}(\zit/p^{r}(n+s)) \label{product}\end{equation}
obtained as the Zariski sheafification of the cup-product map on logarithmic Hodge-Witt cohomology (\ref{cup-prod}) (recall that these sheaves vanish unless $m = n, n+1$ and $q = s, s+1$). There is then an induced cup product map on the higher unramified cohomology groups:
\[ H^{i}_{Zar} (X, \mathcal{H}^{m}_{X}(\zit/p^{r}(n))) \otimes H^{j}_{Zar} (X, \mathcal{H}^{q}_{X}(\zit/p^{r}(s))) \xrightarrow{\cup} H^{i+j}_{Zar} (X, \mathcal{H}^{m+q}_{X}(\zit/p^{r}(n+s)))\]
There is a certain compatibility of the cup product on higher unramified cohomology when $m=n$ with the cup product for logarithmic Hodge-Witt cohomology via the edge maps; i.e., there is a commutative diagram:
\begin{equation}
\begin{tikzcd}
H^{i}_{Zar} (X, \mathcal{H}^{n}_{X}(\zit/p^{r}(n))) \otimes H^{j}_{Zar} (X, \mathcal{H}^{s}_{X}(\zit/p^{r}(s))) \arrow{r}{\cup} \arrow{d}{e^{i,n} \otimes e^{j,s}}  & H^{i+j}_{Zar} (X, \mathcal{H}^{n+s}_{X}(\zit/p^{r}(n+s))) \arrow{d}{e^{i+j,n+s}}\\ H^{i+n}_{\text{\'et}} (X, \zit/p^{r}(n)) \otimes H^{j+s}_{\text{\'et}} (X, \zit/p^{r}(s)) \arrow{r}{\cup}  & H^{i+j+n+s}_{\text{\'et}} (X, \zit/p^{r}(n+s))\label{comm}
\end{tikzcd}
\end{equation}
where the $e^{*,*}$ maps in the vertical arrows are the edge maps. The commutativity of (\ref{comm}) follows from the fact that the Bloch-Ogus spectral sequence has cup products (as in the proof of Proposition \ref{Galois-props}; see also \cite{Sw} Cor. 8.7). 

\subsection{Gersten resolution and consequences}

\noindent Now, let $X^{(i)}$ denote the irreducible closed subsets of $X$ of codimension $i$ and for $x \in X^{(i)}$, let $k(x)$ be the residue field and $i_{x}: \text{Spec } k(x) \hookrightarrow X$ be the corresponding inclusion. Gros and Suwa showed (using Gabber's effacement theorem) that there is an acyclic resolution of $\mathcal{H}^{n}_{X}(\zit/p^{r}(n))$ (\cite{GS} Cor. 1.6) known as the {\em Gersten resolution}:
\begin{equation} 0 \to \mathcal{H}^{n}_{X}(\zit/p^{r}(n)) \to \bigoplus_{x \in X^{(0)}} i_{x*}(W_{r}\Omega^{m}_{log}(k(x))) \xrightarrow{\partial^{0}} \cdots \xrightarrow{\partial^{n-1}}\bigoplus_{x \in X^{(n)}} i_{x*}(\zit/p^{r}) \to 0 \label{Gersten-p}  \end{equation}
where $W_{r}\Omega_{log}^{m}(-)$ is the module of logarithmic Hodge-Witt differentials and the boundary maps
\[ \bigoplus_{x \in X^{(i)}} i_{x*}(W_{r}\Omega^{m-i}_{log}(k(x))) \xrightarrow{\partial^{i}} \bigoplus_{x \in X^{(i+1)}} i_{x*}(W_{r}\Omega^{m-i-1}_{log}(k(x)))\]
are the sums of residue maps over $x \in X^{(i)}$, as explained in the proof of Lem. 4.11 of op. cit. Since (\ref{Gersten-p}) is an acyclic resolution, the Zariski cohomology of $\mathcal{H}^{n}_{X}(\zit/p^{r}(n))$ is given by:
\begin{equation} H^{q}_{Zar} (X, \mathcal{H}^{n}_{X}(\zit/p^{r}(n)))= R^{q}\Gamma ((\ref{Gersten-p}))\label{cohom} \end{equation}
where $\Gamma: Sh(X^{Zar}) \to \text{Ab}$ denotes the global sections functor. Explicitly, $H^{q}_{Zar} (X, \mathcal{H}^{n}_{X}(\zit/p^{r}(n)))$ is the $q^{th}$ cohomology of the complex of Abelian groups:
\begin{equation} C_{X}^{*}(n): \bigoplus_{x \in X^{(0)}} W_{r}\Omega^{m}_{log}(k(x)) \xrightarrow{\partial^{0}} \bigoplus_{x \in X^{(1)}} W_{r}\Omega_{log}^{m-1}(k(x))\xrightarrow{\partial^{1}} \ldots \xrightarrow{\partial^{n-1}} \bigoplus_{x \in X^{(n)}} \zit/p^{r} \label{Gersten-group}  \end{equation}
\noindent A few observations are immediate from (\ref{Gersten-group}). First, we have the vanishing
\begin{equation} H^{q}_{Zar} (X, \mathcal{H}^{n}_{X}(\zit/p^{r}(n))) = 0\label{vanish} \end{equation} for $q>n$ (this is given as Cor. 1.9 of op. cit.). Moreover, there is the so-called Bloch-Quillen formula 
\begin{equation} \psi^{n}: CH^{n} (X)/p^{r} \cong H^{n}_{Zar} (X, \mathcal{H}^{n}_{X}(\zit/p^{r}(n)))\label{BO} \end{equation}
which is Th. 4.13 of op. cit. Using (\ref{BO}), we define the composition with the edge map in the Bloch-Ogus spectral sequence:
\begin{equation} c^{n}: CH^{n} (X)/p^{r} \xrightarrow{\psi^{n}} H^{n}_{Zar} (X, \mathcal{H}^{n}_{X}(\zit/p^{r}(n)))\xrightarrow{e^{n,n}} H^{2n}_{\text{\'et}} (X, \zit/p^{r}(n))\label{cyc-class}\end{equation}
The claim is that $c^{n}$ is the cycle class map in \cite{G} p. 50-51. Indeed, from the proof of \cite{GS} Cor. 1.6, for any irreducible closed subset $x$ of codimension $n$ on $X$ one uses the purity theorem (\cite{G} p. 46) to identify the corresponding $\zit/p^{r}$ summand in (\ref{Gersten-group}) with the relative logarithmic Hodge-Witt cohomology group $H^{n}_{x} (X, \zit/p^{r}(n))$ (note the difference in notation); the map $e^{n,n}$ is that given by the Gysin map, which is precisely how the cycle class map is defined in op. cit. So, in particular, we have compatibility with products:
\begin{equation} c^{m+n} (\alpha\cdot\beta) =  c^{m}(\alpha) \cup c^{n}(\beta) \in H^{2(m+n)}_{\text{\'et}} (X, \zit/p^{r}(m+n)) \label{comp}\end{equation}
where $\alpha$ and $\beta$ are as in the previous paragraph. Additionally, we observe that (\ref{BO}) defines an action of $CH^{*} (X)$ on the higher unramified cohomology groups:
\begin{equation} H^{i}_{Zar} (X, \mathcal{H}^{q}_{X}(\zit/p^{r}(n)) \otimes CH^{m} (X) \xrightarrow{\cdot} H^{i+m}_{Zar} (X, \mathcal{H}^{q+m}_{X}(\zit/p^{r}(n+m)), \ \alpha \cdot \beta := \alpha \cup \psi^{m}(\beta) \label{CH-act} \end{equation}
Finally, we note that when $i=q=n$, the action defined by (\ref{CH-act}) coincides with the usual intersection product for Chow groups (using the identification (\ref{BO})). Indeed, one can use a moving lemma (for instance, \cite{Bl} Lem. 1.1) to move $\beta$ so that its support intersects that of $\alpha$ properly, and then it is easy to see this coincidence.\\
\indent As in \cite{R1} (3.4) and (3.5), one can use the Gersten resolution to define push-forward and pull-back operations. Indeed, for $f: X \to Y$ a proper morphism of relative dimension $d$, there is an induced map of complexes:
\begin{equation} f_{*}: C^{*}_{X} (n) \to C^{*}_{Y} (n-d)\label{naive-push} \end{equation}
obtained in the following way. For $x \in X^{(i)}$, define a map
\begin{equation} \nu_{r}(n-i)(k(x)) \to \nu_{r}(n-i)(k(f(x)))\label{loc-push} \end{equation}
to be the trace map when $[k(x): k(f(x))] < \infty$ and $0$ otherwise. The trace map for field extensions is described in the case of a finite separable extension in the discussion preceding \cite{GS2} Lemma 9.5.4; moreover, an explicit computation of the trace map for a purely inseparable extension is given on pp. 621-623 of \cite{Sh}. Taking the sum over $x \in X^{(i)}$ then gives (\ref{naive-push}). In the case that $i=n$, (\ref{loc-push}) is precisely the map $\zit/p^{r} \to \zit/p^{r}$ given by multiplication by $[k(x): k(f(x))]$ (when it is finite and by $0$ else). Comparing this with the construction of push-forward on Chow groups \cite{F} Ch. 1, we observe that applying $H^{n}_{Zar}$ to (\ref{naive-push}) induces the usual push-forward map on $CH^{n} (-)/p^{r}$, using the identification (\ref{BO}). \\
\indent When $f: X \to Y$ is a flat morphism, there is an induced map of complexes:
\begin{equation} f^{*}: C^{*}_{Y} (n) \to C^{*}_{X} (n) \label{naive-pull}\end{equation}
obtained in the following way. For $y \in Y^{(i)}$ and $x$ be an irreducible component of $f^{-1} (y)$, define the corresponding map
\[\nu_{r}(n-i)(k(y)) \to \nu_{r}(n-i)(k(x))\]
to be the restriction map on differential forms. Taking the sum over $y \in Y^{(i)}$ then gives (\ref{naive-pull}). In the case that $i=n$, (\ref{loc-push}) is precisely the identity map $\zit/p^{r} \to \zit/p^{r}$. As with push-forward, we compare this with the construction of flat pull-back on Chow groups \cite{F} Ch. 1 and observe that applying $H^{n}_{Zar}$ to (\ref{naive-pull}) induces the usual pull-back map on $CH^{n} (-)/p^{r}$, using the identification (\ref{BO}). \\
\indent The operations of flat pull-back and proper push-forward and the action of the Chow group allow one to define the action of correspondences on the higher unramified cohomology groups in the familiar way. With the above definitions of pull-back, pushforward and $\cdot$, for any $\Gamma \in Cor^{q} (X, Y)$ we define a corresponding map
\[ \Gamma_{*}: H^{p}_{Zar} (X, \mathcal{H}^{m}_{X}(\zit/\ell^{r}(n))) \to H^{p+q}_{Zar} (X, \mathcal{H}^{m+q}_{X}(\zit/\ell^{r}(n+q))), \ \Gamma_{*} (\alpha) = \pi_{Y*}(\pi_{X}^{*}\alpha\cdot \psi^{n}(\Gamma)) \]
where $\pi_{X}: X \times Y \to X$ and $\pi_{Y}: X \times Y \to Y$ are the projections. By the compatibilities mentioned above, we have $(\Gamma' \circ \Gamma)_{*} = \Gamma_{*}' \circ \Gamma_{*}$. We would like to check that the action of correspondences defined here is compatible (via the edge map) with action of correspondences on logarithmic Hodge-Witt cohomology:
\begin{Lem}\label{funct} Suppose $X$ and $Y$ be smooth projective schemes over $k$ and that $Y$ has dimension $d_{Y}$. Also, let $\Gamma \in Cor^{q} (X, Y)/p^{r}$. Then, the diagram below commutes:
\[\begin{tikzcd} H^{m}_{Zar} (X, \mathcal{H}_{X}^{n} (\zit/p^{r}(n)))) \arrow{d}{e^{m,n}} \arrow{r}{\Gamma_{*}} & H^{m+q}_{Zar} (Y, \mathcal{H}_{Y}^{n+q} (\zit/p^{r}(n+q)))) \arrow{d}{e^{m+q,n+q}} \\ 
H^{m+n}_{\text{\'et}} (X, \zit/p^{r}(n))) \arrow{r}{c^{q+d_{Y}}(\Gamma)_{*}}& H^{m+n+2q}_{\text{\'et}} (Y, \zit/p^{r}(n+q)))\end{tikzcd} \]
for $m =n-1,n$.
\begin{proof} By the definition of $\Gamma_{*}$, it will suffice to verify that $e^{*,*}$ commutes with proper push-forward, flat pull-back and products. That the edge map commutes with products was already noted in (\ref{comm}). The commutativity of the edge map with push-forward and pull-back is addressed in \cite{B} Lemma 3.2 for \'etale cohomology with $\zit/\ell$ coefficients (when $\ell$ is different from the characteristic); the commutativity is a formal consequence of functoriality/duality properties (enumerated in \cite{BO}) for this cohomology theory. The analogous functoriality/duality properties for logarithmic Hodge-Witt cohomology are established by Gros in \cite{G}. 
\end{proof}
\end{Lem}

\subsection{A Bockstein exact sequence} 
\noindent Let $CH^{n} (-, q)$ denote Bloch's higher Chow groups \cite{Bl}. For $X$ a smooth quasi-projective scheme over $k$, we let $\C\H^{n}_{X}(q)$ denote the corresponding Zariski sheafification of the presheaf of Abelian groups over $X$ defined by $U \mapsto CH^{n} (U, q)$.
We note that there is a Gersten resolution by op. cit. \S 10:
\begin{equation} 0 \to \C\H^{n}_{X}(q)  \to \bigoplus_{x \in X^{(0)}} i_{x*}(CH^{n}(k(x), q)) \xrightarrow{\partial^{0}} \cdots \xrightarrow{\partial^{n-1}}\bigoplus_{x \in X^{(n)}} i_{x*}(CH^{0}(k(x), q-n)) \to 0 \label{Gersten-B}  \end{equation}
In particular, it follows that $\C\H^{n}_{X}(q) = 0$ for $n>q$ using the fact that $CH^{n}(k(x), q) = 0$ for $n>q$ (for degree reasons). Also, when $n=q$, the (independent) results of Nesternko-Suslin \cite{NS} and Totaro \cite{T} then show that there is a natural isomorphism $K_{n}^{M} (-) \cong CH^{n} (-, n)$, where $K_{n}^{M}$ denotes the $n^{th}$ Milnor $K$-theory \cite{Mi}. For ease of notation, we let $\C\H^{n,n}_{X} := \C\H^{n}_{X}(n)$. When $k$ is of characteristic $p>0$, Geisser and Levine establish the short exact sequence of Zariski sheaves below (see \cite{GL} proof of Theorem 8.5):
\begin{equation} 0 \to \C\H^{n,n}_{X} \xrightarrow{\times p^{r}} \C\H^{n,n}_{X} \to \alpha_{*}\nu_{r}(n) = \mathcal{H}^{n}_{X}(\zit/p^{r}(n))\to 0\label{GL-short} \end{equation}
\begin{Prop}\label{big} Let $X$ be a smooth quasi-projective scheme over $k$. Then, there exists a natural short exact sequence:
\begin{equation} 0 \to H^{n-1}_{Zar} (X, \C\H^{n,n}_{X})/p^{r} \to H^{n-1}_{Zar} (X, \mathcal{H}^{n}_{X}(\zit/p^{r}(n))) \xrightarrow{\delta_{X}^{n}} CH^{n} (X)[p^{r}] \to 0\label{short}\end{equation}
for which $\delta_{X}^{n}$ commutes with flat pull-back and proper push-forward and is compatible with the action of the Chow group.
\begin{proof} By applying $H^{*}_{Zar} (X, -)$ to (\ref{GL-short}) and truncating, one obtains the short exact sequence:
\[ 0 \to H^{n-1}_{Zar} (X, \C\H^{n,n}_{X})/p^{r} \to H^{n-1}_{Zar} (X, \mathcal{H}^{n}_{X}(\zit/p^{r}(n))) \xrightarrow{\delta_{X}^{n}} H^{n}_{Zar} (X, \C\H^{n,n}_{X})[p^{r}] \to 0\]
The identification $CH^{n} (X) \cong H^{n}_{Zar} (X,  \C\H^{n,n}_{X})$ follows from the corollary on p. 300 of \cite{Bl}. The functoriality statements then follow from the corresponding statements for (\ref{GL-short}) or rather for the map
\begin{equation} \C\H^{n,n}_{X} \to \alpha_{*}\nu_{r}(n)\label{bound} \end{equation}
in (\ref{GL-short}). By the Gersten resolution, we immediately reduce to the case of a field $F$ of characteristic $p>0$; i.e., we would like to show that the map
\begin{equation} CH^{n} (F, n) \xleftarrow{\cong} K_{n}^{M} (F) \xrightarrow{dlog_{n}} W_{r}\Omega^{n}_{log}(F) \label{log-map}\end{equation}
defining (\ref{bound}) is compatible with restriction maps, is multiplicative and carries the norm on the left to the trace on the right. To establish this, we first describe the map (\ref{log-map}). For $x_{1}, \ldots x_{n} \in k^{\times}$ we denote the corresponding symbol by
\[ \{ x_{1}, \ldots x_{n} \} \in  K_{n}^{m}(F) \]
The {\em log-symbol map} which we denote by $dlog_{n}$ is defined by $\{ x_{1}, \ldots x_{n} \} \mapsto \frac{dx_{1}}{x_{1}} \wedge \ldots \wedge \frac{dx_{n}}{x_{n}}$ when $r=1$ (and there is the analogous definition involving Teichm\"uller lifts of the $x_{i}$ when $r>1$). There are versions of the isomorphism in (\ref{log-map}) that depend on whether one uses the simplicial \cite{NS} or the cubical \cite{T} definition of Bloch's cycle complex. We will use the cubical version, which is defined as 
\begin{equation} \{ x_{1}, \ldots x_{n} \} \mapsto (\frac{x_{1}}{x_{1}-1}, \ldots, \frac{x_{n}}{x_{n}-1}) \label{Tot}  \end{equation}
where the right hand side is viewed as the class in $CH^{n} (F, n)$ of an $n$-tuple of $(\mathbb{P}^{1}\setminus 1)^{n}_{F}$. That (\ref{Tot}) and $dlog_{n}$ are multiplicative and compatible with restriction is evident from the definition. That $dlog_{n}$ has the desired property with respect to norm and trace maps was established by Kato in \cite{Ka} \S 3.3 Lemma 13. Finally, that (\ref{Tot}) is compatible with norm follows easily from Totaro's definition of the inverse of (\ref{Tot}) in \cite{T} p. 183. Indeed, this definition involves the norm map in $K$-theory and using the functoriality property mentioned in \cite{BT} p. 385, the desired compatibility follows.
\end{proof}
\end{Prop}
\begin{Cor}\label{cor-big} For $X$ a smooth projective scheme over a perfect field $k$ of characteristic $p>0$, the boundary map $H^{n-1}_{Zar} (X, \mathcal{H}^{n}_{X}(\zit/p^{r}(n))) \xrightarrow{\delta_{X}^{n}} CH^{n} (X)[p^{r}]$ commutes with the action of correspondences.
\end{Cor}
\begin{Rem} The well-known Bloch-Gabber-Kato Theorem \cite{BK} shows that the log-symbol is in fact an isomorphism after modding out by $p^{r}$. However, we have chosen to state the above result using Bloch's higher Chow groups because of the Gersten resolution, which is harder to establish for Milnor $K$-theory (although Kerz has shown that this holds over infinite fields; see \cite{Ke}).
\end{Rem}

\section{Rost nilpotence principle}
\subsection{Main result}
\begin{Not} For $E/k$ a field extension, consider the extension-of-scalars map:
\[ CH^{i} (X) \to CH^{i} (X_{E}) \]
For any $\alpha \in CH^{i} (X)$, we will denote by $\alpha_{E}$ its image under this map.
\end{Not}

\begin{Conj}[Rost nilpotence principle] Let $X$ be a smooth projective scheme over a field $k$ and $\gamma \in Cor^{0} (X, X)$. If $\gamma_{F} = 0$ for some field extension $F/k$, then $\gamma^{\circ N} = 0$ for $N \gg 0$.
\end{Conj}

\noindent To prove this conjecture in the case of surfaces, we begin with a straightforward observation. Since the proof is standard, we only sketch it.
\begin{Lem}\label{above} Let $k$ be perfect, $F/k$ a field extension, $X$ a smooth projective scheme over $k$ and $\alpha \in CH^{i} (X)$. If $\alpha_{F} = 0$, then $\alpha$ is torsion and $\alpha_{\overline{k}} = 0$.
\begin{proof} To prove that $\alpha$ is torsion, we can take $F$ to be a finitely generated over $k$. Then, by \cite{Hart} Theorem I.4.8A, $F$ is separably generated; i.e., there is some finitely generated purely transcendental extension $K/k$ for which $F/K$ is finite and separable. 
By a standard transfer argument (see, for instance, \cite{F} Ex. 1.7.4), we reduce to the case that $K$ is purely transcendental over $k$. As noted by one of the referees, a reference for the well-known fact that Chow groups are preserved under purely transcendental extensions can be found in \cite{FOV} Proposition 2.1.8. Hence, $\alpha$ is torsion. To prove that $\alpha_{\overline{k}} = 0$, let $F/k$ be a field extension and $\overline{F}$ the algebraic closure of $F$. Then, consider the commutative diagram below:
\[\begin{tikzcd}
CH^{i} (X)_{tors} \arrow{r} \arrow{d} & CH^{i} (X_{F}) \arrow{d}\\
CH^{i} (\overline{X})_{tors} \arrow{r} & CH^{i} (X_{\overline{F}}) 
\end{tikzcd}
\]
By Theorem 3.11 of \cite{Le}, the bottom horizontal arrow is injective. Thus, the kernel of the left vertical arrow contains the kernel of the top horizontal arrow. So, if $\alpha_{F} = 0$, then $\alpha_{\overline{k}} = 0$, as desired.
\end{proof}
\end{Lem}
\noindent Lemma \ref{above} implies that if $\gamma \in Cor^{0} (X, X)$ satisfies the hypothesis of the Rost nilpotence principle, we can assume without loss of generality that $F = \overline{k}$ and that $\gamma$ is torsion. Since
\[ CH^{d} (X\times X)_{tors} = \bigoplus_{\ell \text{ prime}} CH^{d} (X\times X)[\ell^{\infty}] \]
we can write $CH^{d} (X\times X)_{tors} = CH^{d} (X\times X)' \oplus CH^{d} (X\times X)[p^{\infty}]$. Thus, we will focus on the case that $\gamma$ is a $p$-primary torsion cycle, since the case of prime-to-$p$ torsion cycles have already been handled in \cite{RS}.
\begin{Thm}\label{main} Let $X$ be a smooth projective scheme of dimension $d$ over a perfect field $k$ of characteristic $p>0$ and $\gamma \in CH^{d} (X\times X)[p^{\infty}]$. Suppose that $\gamma_{\overline{k}} = 0$ and that
\begin{equation} H^{d-3}_{Zar} (X \times X, \mathcal{H}^{d+1}_{X \times X}(\zit/p^{r}(d))) = 0\label{obstr} \end{equation}
Then, $\gamma^{n} = 0 \in Cor^{0} (X, X)$ for $n \gg 0$.
\begin{proof}
For convenience set $Y = X \times X$ and let $\tilde{\gamma} \in H^{d-1}_{Zar} (Y, \mathcal{H}^{d}_{Y}(\zit/p^{r}(d)))$ be a lift of $\gamma$ via the boundary map:
\[ \begin{tikzcd}\delta_{Y}^{d}: H^{d-1}_{Zar} (Y, \mathcal{H}^{d}_{Y}(\zit/p^{r}(d))) \arrow[two heads]{r} & CH^{d} (Y)[p^{r}]\end{tikzcd} \]
Then, consider the correspondence $\beta_{N} = \gamma^{\circ N} \times \Delta_{Y} \in Cor^{0} (Y, Y).$
We note that $\beta_{N*} (\gamma) = \gamma \circ \gamma^{\circ N}  = \gamma^{\circ N+1}$ using \cite{MNP} Lem. 2.1.3. Now, it follows from Corollary \ref{cor-big} that
\[ \delta_{Y}^{d} (\beta_{N*} (\tilde{\gamma})) = \beta_{N*} (\delta_{Y}^{d} (\tilde{\gamma})) = \gamma^{\circ N+1} \in  CH^{d} (Y)[p^{r}]\] 
To show that this is $0$, it suffices to show that 
\begin{equation} \beta_{N*}: H^{d-1}_{Zar} (Y, \mathcal{H}^{d}_{Y}(\zit/p^{r}(d))) \to H^{d-1}_{Zar} (Y, \mathcal{H}^{d}_{Y}(\zit/p^{r}(d)))\label{target} \end{equation}
vanishes for $N \gg 0$. Indeed, note that 
\begin{equation} c^{d} (\gamma)^{\circ N} = c^{d} (\gamma^{\circ N}) = 0 \in H^{2d}_{\text{\'et}} (Y, \zit/p^{r}(d))\label{beta} \end{equation}
for $N\geq 2d+1$ by Corollary \ref{cor-van}, from which we deduce that
\begin{equation} \beta_{N*}: H^{*}_{\text{\'et}} (Y, \zit/p^{r}(d)) \to H^{*}_{\text{\'et}} (Y, \zit/p^{r}(d))\label{betaN}\end{equation}
vanishes for $N \gg 0$. To prove the result, note that by the long exact sequence (\ref{long}), there is a short exact sequence:
\begin{equation} H^{d-3}_{Zar} (Y, \mathcal{H}^{d+1}_{Y}(\zit/p^{r}(d))) \to H^{d-1}_{Zar} (Y, \mathcal{H}^{d}_{Y}(\zit/p^{r}(d))) \xrightarrow{e^{d-1,d}} H^{2d-1}_{\text{\'et}} (Y, \zit/p^{r}(d)) \label{last-ses}\end{equation}
Since $e^{d-1,d}$ commutes with correspondences by Lemma \ref{funct}, $\beta_{N}$ acts on (\ref{last-ses}). Since the first term in (\ref{last-ses}) vanishes by assumption, it follows from (\ref{betaN}) that (\ref{target}) vanishes for $N \gg 0$, as desired.
\end{proof}
\end{Thm}
\noindent In the case that $d=2$, we obtain the following consequence.
\begin{Cor}\label{surf} The Rost nilpotence principle holds for a smooth projective surface $S$ over a perfect field.
\begin{proof} 
As mentioned previously, we reduce to the case that $\gamma_{\overline{k}} = 0$ and that $\gamma$ is torsion. The case that $\gamma$ is prime-to-p torsion case was already proved in \cite{RS}. Indeed, Theorem 1.1 of op. cit. shows that $\gamma^{\circ n} = 0 \in CH^{2}_{L} (S \times S)$ for $n \gg 0$, where $CH^{*}_{L}$ denotes the Lichtenbaum Chow group (see op. cit. for a review of this notion). Since the natural map $CH^{2} (S \times S) \to CH^{2}_{L} (S \times S)$ is injective by \cite{K} Prop. 2.9, it follows that $\gamma^{\circ n} = 0$. When $\gamma$ is a $p$-primary torsion cycle, it follows immediately from Theorem \ref{main} that $\gamma^{\circ n} = 0$ since (\ref{obstr}) vanishes trivially in the case that $d=2$.
\end{proof}
\end{Cor}

\subsection{Decompositions of the diagonal}

\begin{Def}\label{def} Given a scheme $X$ over a field $k$, we say that the Chow group of $X$ is {\em universally supported in dimension $\leq i$} if there exists a Zariski closed subset $V \subset X$, all of whose irreducible components are of dimension $\leq i$, for which the push-forward
\[ CH_{0} (V_{F}) \to CH_{0} (X_{F}) \]
is surjective for all field extensions $F/k$. 
\end{Def}
\noindent The Bloch-Srinivas decomposition method (as in the proof of \cite{BS} Prop. 1) gives the following result directly:
\begin{Lem}\label{gamma-lem} Suppose that $X$ is a smooth projective scheme of dimension $d$ over a field $k$ whose Chow group is {universally supported in dimension $\leq i$}. Then, there exists a divisor $D \subset X$ and a closed subset $V \subset X$ (as in the above definition) such that for every  $\gamma \in CH^{d} (X \times X)$
\begin{equation} \gamma =  \gamma_{1} + \gamma_{2} \in CH^{d} (X \times X)\label{gamma-decomp}\end{equation}
for some $\gamma_{1}$ is supported on $D \times X$ and $\gamma_{2}$ is supported on $X \times V$. Moreover, if for some field extension $F/k$ we have 
\[\gamma_{F} = 0 \in CH^{d} (X_{F} \times X_{F})\] 
then one can take $\gamma_{i}$ to be such that $\gamma_{i, F} = 0$ for $i=1,2$.
\begin{proof} Only the second statement requires justification, given op. cit. For this, we observe by \cite{MNP} Lem. 2.1.3 that
\[  \gamma = (\Delta_{X} \times \gamma)_{*}(\Delta_{X}) \in CH^{d} (X \times X) \]
Then, $\Delta_{X}$ decomposes as $\Delta_{1} +\Delta_{2}$ as in (\ref{gamma-decomp}) so that we can take $\gamma_{i} = (\Delta_{X} \times \gamma)_{*}(\Delta_{i})$, for which we have $\gamma_{i, F} = 0$.
\end{proof}
\end{Lem}
\noindent In the lemma below, we will need to assume the following:
\begin{Ass*}[Resolution of singularities in dimension $\leq d$] Let $X$ be a projective scheme over a perfect field $k$, all of whose irreducible components are of dimension $\leq d$. Then, there exists a smooth projective scheme $\tilde{X}$ over $k$ and a surjective morphism $\phi: \tilde{X} \to X$ that is birational over each irreducible component of $X$ (and is an isomorphism over the non-singular locus of each component).
\end{Ass*}
\noindent Of course, this assumption holds for all $d$ when $k$ is a field of characteristic zero by Hironaka's well-known theorem \cite{Hi}. In positive characteristic, it is known to hold when $d \leq 3$ by the main result of \cite{CP} (building on classical work of Zariski, Abhyankar and others). We are now able to state the following lemma that allows one to (partially) bootstrap the Rost nilpotence principle to higher dimension.
\begin{Lem}\label{supp} Assume that the following hold:
\begin{enumerate}[label=(\alph*)]
\item\label{Rost} the Rost nilpotence principle for irreducible smooth projective schemes of dimension $\leq d-1$;
\item resolution of singularities in dimension $\leq d-1$.
\end{enumerate}
Let $X$ be a smooth projective scheme of dimension $d$ whose Chow group is universally supported in dimension $\leq d-1$. Then, $X$ satisfies the Rost nilpotence principle.
\begin{proof} Suppose that $\gamma \in Cor^{0} (X, X)$ is a correspondence for which $\gamma_{F} = 0$. By Lemma \ref{gamma-lem}, there is a decomposition (\ref{gamma-decomp}) for which $\gamma_{i, F} = 0$. Since a sum of nilpotent correspondences is again nilpotent, we reduce to showing that each $\gamma_{i}$ is nilpotent. By the symmetry of the problem, it suffices to prove the nilpotency of $\gamma_{2}$; i.e., we may assume that $\gamma$ lies in the image of the pushforward:
\[ CH^{d-1} (X \times V) \to CH^{d} (X \times X) \]
where $V \subset X$ is a closed subset whose irreducible components all have dimension $\leq d-1$. Now, because of the resolution of singularities assumption, there exists a smooth projective scheme $\tilde{V}$ over $k$ and a morphism $\tilde{V} \to V$ for which the conclusion of Lemma \ref{anc-2} holds. In particular, letting $\phi: \tilde{V} \to V \hookrightarrow X$ be the composition, we can write
\[ \gamma = (\Delta_{X} \times \phi)_{*}(\alpha) = \Gamma_{\phi} \circ \alpha \in CH^{d} (X \times X)\]
for $\alpha \in CH^{d-1} (X \times \tilde{V})$, again using \cite{MNP} Lemma 2.1.3. Now, we consider 
\[ \beta= \alpha \circ  \Gamma_{\phi} \in Cor^{0} (\tilde{V}, \tilde{V})\]
Then, we have $\beta^{\circ 2} = \alpha \circ \gamma \circ \Gamma_{\phi}$. Since $\gamma_{F} = 0$, it follows that $\beta^{\circ 2}_{F} = 0$. By Lemma \ref{anc-1} (and assumption \ref{Rost}), the Rost nilpotence principle holds for $\tilde{V}$. We deduce that $\beta^{\circ 2}$ is nilpotent and, hence, so is $\beta$. Since for $n \geq 0$ we have
\[ \gamma^{\circ n+1} = \Gamma_{\phi} \circ \beta^{\circ n} \circ \alpha  \]
it follows that $\gamma$ is also nilpotent, as desired.
\end{proof}
\end{Lem}
\begin{Rem} The proof of the above lemma may easily be modified to an induction argument for the Rost nilpotence principle in general, provided that one is able to prove that for $F/k$ a field extension and
\[ \gamma \in \text{Ker }\{ Cor^{0} (X, X) \to Cor^{0} (X_{F} \times X_{F}) \} \] 
some power of $\gamma$ admits a decomposition such as (\ref{gamma-decomp}). 
\end{Rem}
\noindent We obtain the following as an immediate consequence of Lemma \ref{supp}, Corollary \ref{surf} and the fact that resolution of singularities is known in dimension $\leq 2$ (by the main results of \cite{Ab}, \cite{Li}):
\begin{Cor}\label{three} Suppose that $X$ is a smooth projective scheme of dimension $\leq 3$ over a perfect field $k$ whose Chow group is universally supported in dimension $\leq 2$. Then, $X$ satisfies the Rost nilpotence principle.
\end{Cor}
\noindent Recall that a smooth projective scheme $X$ of dimension $d$ over a field $k$ is said to be {\em birationally ruled} if it is irreducible and is $k$-birational to $Y \times \pit^{1}$, where $Y$ is a smooth projective scheme of dimension $d-1$. Rosenschon and Sawant in \cite{RS} prove the Rost nilpotence principle for birationally ruled threefolds in characteristic $0$. Their proof does not generalize so easily to positive characteristic, since it invokes a non-trivial result in birational geometry known as the Weak Factorization Theorem. As an application of Corollary \ref{three}, we can prove the following generalization of their result:
\begin{Cor} Suppose that $X$ is a birationally ruled threefold over a perfect field. Then, $X$ satisfies the Rost nilpotence principle.
\begin{proof} By the localization sequence \cite{F} \S 1.8, it suffices to prove that there exists a Zariski open subset $U \subset X$ for which $CH_{0} (U_{F}) =0$ for all field extensions $F/k$. By assumption, there exists a some smooth projective surface $S$ and a rational map (defined over $k$) $\phi: X \dashrightarrow S \times \pit^{1}$ which is an isomorphism $U \xrightarrow[\cong]{\phi} V$ for some open subsets $U \subset X$ and $V \subset S \times \pit^{1}$. Now, by the projective bundle formula, the Chow group of $S \times \pit^{1}$ is universally supported in dimension $\leq 2$. So, by shrinking $V$ if necessary, we deduce from the localization sequence that $CH_{0} (V_{F}) = 0$ for all field extensions $F/k$. It follows that $CH_{0} (U_{F}) =0$, as desired.
\end{proof}
\end{Cor}

\begin{appendices}
\section{Two ancillary lemmas} 

\noindent Here we state and prove two routine lemmas used in the previous section. The first concerns the Rost nilpotence principle in the case of a smooth projective scheme that is not equi-dimensional. 
\begin{Lem}\label{anc-1} Suppose that the Rost nilpotence principle holds for all integral smooth projective schemes of dimension $\leq d$ over a field $k$. Then, it also holds for smooth projective schemes over $k$ whose irreducible components are all of dimension $\leq d$.
\begin{proof} Let $X$ be a smooth projective scheme over $k$ whose irreducible components are all of dimension $\leq d$. We induct on the number of irreducible components $n$ of $X$. The base case that $n=1$ is true by assumption. So, by way of induction, we assume that Rost nilpotence holds for all smooth projective schemes over $k$ with $n-1$ irreducible components (all of dimension $\leq d$). If $X$ has exactly $n$ irreducible components, write $X = X_{0} \coprod X_{1}$, where $X_{0}$ has $n-1$ irreducible components and $X_{1}$ is irreducible. Then, we have
\[ Cor^{0} (X, X) = Cor^{0} (X_{0}, X_{0}) \oplus Cor^{0} (X_{1}, X_{0}) \oplus Cor^{0} (X_{0}, X_{1}) \oplus Cor^{0} (X_{1}, X_{1})  \]
Suppose that $\gamma \in Cor^{0} (X, X)$ is such that $\gamma_{F} = 0$ for some field extension $F/k$. Then, writing
\[ \gamma = \gamma_{00} + \gamma_{10} + \gamma_{01} + \gamma_{11} \]
where $\gamma_{ij} \in Cor^{0} (X_{i}, X_{j})$, we have $\gamma_{ij, F} = 0$. The inductive hypothesis then implies that
\begin{equation} \gamma_{ii}^{\circ n} = 0 \in Cor^{0}(X_{i}, X_{i}) \label{ind}\end{equation}
Moreover, by definition of $\circ$, we have $\gamma_{kl} \circ \gamma_{ij} = 0$ unless $j=k$. Thus, we have
\begin{equation} \gamma^{\circ 2} = \sum_{i,j,k} \gamma_{jk}\circ\gamma_{ij}\label{sum} \end{equation}
where each of the indices $i,j,k$ is either $0$ or $1$. Since there are only two possible indices, we must have either $i=j$ or $j=k$ or $i=k$ in all of the summands above. In the case that $i=j$, the summand $\gamma_{ik}\circ\gamma_{ii}$ above is nilpotent; indeed, $\gamma_{ii}$ is nilpotent by (\ref{ind}), which means that $\gamma_{ik}\circ\gamma_{ii}$ is also (the set of nilpotent elements forms a two-sided ideal). The same is true if $j=k$. On the other hand, if we have $i =k$, then the summand $\gamma_{ji}\circ\gamma_{ij}$ vanishes upon extending scalars to $F$, since each of $\gamma_{ij}$ and $\gamma_{ji}$ does. This means that 
\[ \gamma_{ji}\circ\gamma_{ij} \in Cor^{0}(X_{i}, X_{i})\]
is nilpotent by the inductive hypothesis. So, all the summands on the right hand side (\ref{sum}) are nilpotent, which means that so is $\gamma^{\circ 2}$.
\end{proof}
\end{Lem}
\noindent There is also the following Chow-theoretic consequence of the resolution of singularities assumption. The proof is a standard inductive argument involving the localization sequence, and its proof has been left to the reader (in fact, as one of the referees mentioned, this type of argument is given in \cite{Voi} Theorem 9.27).
\begin{Lem}\label{anc-2} Suppose that $k$ is a perfect field for which resolution of singularities in dimension $\leq d$ holds. Let $V$ be a projective scheme over $k$ all of whose irreducible components have dimension $\leq d$. Then, there exists some smooth projective scheme $V'$ over $k$ all of whose irreducible components have dimension $\leq d$ and a morphism $\phi: V' \to V$ for which the push-forward
\[ (\text{id}_{X} \times \phi)_{*}: CH_{*} (X \times V') \to CH_{*} (X \times V) \]
is surjective for all schemes $X$ over $k$. 
\end{Lem}
\end{appendices}

\Addresses
\end{document}